\documentclass[12pt]{article}

\usepackage{fullpage}
\usepackage[usenames]{color}
\usepackage[colorlinks=true,
linkcolor=webgreen,
filecolor=webbrown,
citecolor=webgreen]{hyperref}

\definecolor{webgreen}{rgb}{0,.5,0}
\definecolor{webbrown}{rgb}{.6,0,0}
\definecolor{red}{rgb}{1,0,0}

\usepackage{amsthm, amsfonts, amsmath, amssymb}
\usepackage{mathtools}
\usepackage[margin=1.00in]{geometry}
\usepackage{hyphenat}
\usepackage{float}
\usepackage{graphicx}
\usepackage{tabularx}

\newtheorem{thm}{Theorem}[section]

\newtheorem{defn}{Definition}[section]
\newtheorem{lem}{Lemma}[section]
\newtheorem{conj}{Conjecture}[section]
\newtheorem{coro}{Corollary}[section]
\newtheorem{remark}{Remark}[section]

\def\modd#1 #2{#1\ \mbox{\rm (mod}\ #2\mbox{\rm )}}

\usepackage{hyperref}
\newcommand{\seqnum}[1]{\href{http://oeis.org/#1}{\underline{#1}}}

\usepackage{currfile}
\usepackage{url}
\usepackage{enumitem}
\usepackage{sectsty}

\sectionfont{\fontsize{11}{15}\selectfont}

\allowdisplaybreaks

\begin{document}

\begin{center}
\vskip 1cm{\LARGE\bf 
A Collatz Conjecture Proof \\
}
\vskip 1cm
\large
Robert H. Nichols, Jr.\\ 
The Labry School of Science, Technology, and Business\\
Cumberland University \\ 
Lebanon, TN 37087\\
USA \\
\href{mailto:rnichols@cumberland.edu}{\tt rnichols@cumberland.edu} \\
\today
\end{center}

\vskip .2 in

\begin{abstract}
We represent the generalized Collatz function with the recursive ruler function $r(2n)=r(n)+1$ and $r(2n+1)=1$. We generate even-only and odd-only Collatz subsequences that contain significantly fewer elements term by term, to $2$ and $1$, respectively, than are present in the original $3n+1$ and the Terras-modified Collatz sequences. We show that a nonlinear, coupled system of difference equations yields a complete acyclic (except for the trivial cycle) Collatz tree in odds not divisible by $3$ with root vertex $1$. We construct a complete Collatz tree with the axiom of choice and prove the Collatz conjecture.
\end{abstract}

\section{Introduction: The Collatz sequence and the Collatz conjecture}
We address a long--standing and still--open mathematical problem regarding integer sequences attibuted to Lothar Collatz in the 1930s (cf. \cite{Lag2011}).  We represent the non-negative integers with $\mathbb{N}$ and the positive integers with $\mathbb{N}_{>0}$. We define a $sequence$ and a $subsequence$ of positive integers.
\begin{defn}[Sequence]\label{seq}
For $n \in \mathbb{N}$ and $x_n \in \mathbb{N}_{>0}$ a sequence  $X$ in $\mathbb{N}_{>0}$ is an ordered set of positive integers denoted by X = $X(x_0) = (x_n )_{n \ge 0} = (x_0, x_1, \ldots)$.
\end{defn}
\begin{defn}[Subsequence]\label{subseq}
For $n \in \mathbb{N}$ and  $x_n \in \mathbb{N}_{>0}$  let $(x_n )_{n \ge 0}$ be a sequence. A subsequence of $(x_n )_{n \ge 0}$ is a sequence $( x_{n_i} )_{i \ge 0}$ with $n_i < n_{i+1}$ for all $i \ge 0$.
\end{defn}
If a sequence $A$ is a subsequence of $B$, then we say that $B$ \textit{covers} $A$ and denote the \textit{covering relation} by $A \lessdot B$. For a sequence $A$ the \textit{cardinality} $|A|$ is the number of elements in $A$. Definitions~\ref{seq} and \ref{subseq} allow sequences, for example, to have infinite cardinality, to be non-decreasing, or to contain \textit{cycles} of repeating elements.  If we choose to terminate a sequence $A$ at the first occurrence of a particular number, say $f \in \mathbb{N}$, then we denote the \textit{finite sequence} as $A_f$ with cardinality $|A_f|$.  The smallest $n$ such that $a_n = f$ is the \textit{stopping time} of $A(a_0)$ and is denoted by $A_S$. Thus, if we have a finite sequence $A_f$ then we have that $|A_f(a_0)| = A_S(a_0) + 1$.

For $n \in \mathbb{N}$ and $c_n \in \mathbb{N}_{>0}$  we define the \textit{Collatz sequence} $C (c_0)= (c_0, c_1, \ldots)$ with the recursive function
\begin{equation}\label{Cseq}
c_{n+1} =
\begin{cases}
c_n/2 & \text{if } c_n  \equiv \modd{0} {2} \\
3c_n+1 &  \text{if } c_n \equiv \modd{1} {2} \\
\end{cases}
\end{equation}
(cf. \cite{Lag2011, Lag2012}, \seqnum{A014682} \cite{Sloane2018}, \seqnum{A070165} \cite{Sloane2018}).  The \textit{Collatz conjecture} (\textit{CC}) is that $1 \in C(n)$ for all $n \in \mathbb{N}_{>0}$ (cf. \cite{Lag2011, Lag2012}). Other formulations of the \textit{CC} exist (cf. \cite{Lag2011, Lag2012}). The \textit{CC} is equivalent to the statements that $|C_1(n)| $ is finite for all $n \in \mathbb{N}_{>0}$ \cite{Terras1976} and that the repeating subsequence $(4,2,1)$ is the only cycle in $C$.  In addition,  if ``sufficient'' subsequences of the integers satisfy the \textit{CC} then the \textit{CC} is true (see \cite{Lag2011, Lag2012}, e.g. \cite{Monks2006}).

\section{The Terras-modified Collatz sequence}

Terras (1976) \cite{Terras1976} defines a subsequence of $C$ \eqref{Cseq} to address the Collatz conjecture. For $n \in \mathbb{N}$ and $t_n \in \mathbb{N}_{>0}$ the \textit{Terras-modified Collatz sequence} $T(t_0) = (t_0, t_1, \ldots)$ is
\begin{equation}\label{Tseq}
t_{n+1} =
\begin{cases}
t_n/2 & \text{if } t_n  \equiv \modd{0} {2} \\
(3t_n+1)/2 &  \text{if } t_n \equiv \modd{1} {2} \\
\end{cases}
\end{equation}
 (\cite{Terras1976, Terras1979}, cf. \cite{Lag2011, Lag2012}, \seqnum{A070168} \cite{Sloane2018}).  For $n \in \mathbb{N}_{>0}$ we have by construction that $T(n)$~\eqref{Tseq} is a subsequence of $C(n)$~\eqref{Cseq}.  Further, by construction as well, we have that if sequences $T$ satisfy the Collatz conjecture, then the Collatz conjecture is true and vice versa. 

\section{A generalized alternating parity Terras-modified Collatz sequence}

For $n \in \mathbb{N}$, $g_n \in \mathbb{N}_{>0}$,  odd $a \in \mathbb{N}$, and odd $b \in \mathbb{N}$ we define a \textit{generalized Terras-modified Collatz sequence} $G(g_0) = (g_0, g_1, \ldots)$ with the recursive function
\begin{equation}\label{Gseq}
g_{n+1} =
\begin{cases}
g_n/2 & \text{if } g_n  \equiv \modd{0} {2} \\
(a g_n+b)/2 &  \text{if } g_n \equiv \modd{1} {2}. \\
\end{cases}
\end{equation}
We note that Conway (\cite{Conway1972}, cf. \cite{Lag2011, Conway2013}) has addressed some similar generalizations of the Collatz problem.

We now solve the difference equations that generate sequences $G$ \eqref{Gseq}. We first define the $parity$ of an integer as its property of being $even$ or $odd$. We present these results in the following theorem.
\begin{thm}\label{thmRecursion}
Let $\mu(g_0) \in \mathbb{N}$ and $\nu(g_0) \in \mathbb{N}$ be the number of iterations of \eqref{Gseq} for even and odd $g_0 \in G$ to change parity to odd $g_{\mu(g_0)}$ and even $g_{\nu(g_0)}$, respectively. For even $g_0$ the odd solution to \eqref{Gseq} is
\begin{equation}\label{evenparity}
g_{\mu(g_0)} =
\bigg(\dfrac{1}{2}\bigg)^{\mu(g_0)}g_0.
\end{equation}
For odd $g_0$ the even solution to \eqref{Gseq} is
\begin{equation}\label{oddparity}
g_{\nu(g_0)} =
\bigg(\dfrac{a}{2}\bigg)^{\nu(g_0)}
\bigg(g_0 + \dfrac{b}{a-2}\bigg) - \dfrac{b}{a-2}
\end{equation}
with $b \equiv 0$ \textnormal{mod ($a-2$)}.
\begin{proof}
First, we consider the initial condition that $g_0 \equiv 0$ (mod $2$) in $\eqref{Gseq}$. The associated first-order, linear, homogeneous difference equation is
\begin{equation}\label{evendiffeqn}
2g_{n+1} - g_n = 0.
\end{equation}
We assume a solution to \eqref{evendiffeqn} of the form $g_n = m\lambda^n$ with $m \in \mathbb{N}$, $n \in \mathbb{N}$, and $\lambda \in \mathbb{R}$. We find that $m = g_0$ and $\lambda = 1/2$, yielding
\begin{equation}\label{evendiffeqnsoln}
g_n = \bigg(\dfrac{1}{2}\bigg)^n g_0.
\end{equation}
For even $g_0$ in \eqref{evendiffeqnsoln} we have that $g_n \equiv 1$ (mod $2$) when $n=\mu(g_0)$, yielding \eqref{evenparity}.  Second, we consider the initial condition that $g_0 \equiv 1$ (mod $2$) in \eqref{Gseq}.  The associated first-order, linear, inhomogeneous difference equation is
\begin{equation}\label{odddiffeqn}
2g_{n+1} - a g_n = b.
\end{equation}
To find homogeneous and particular solutions to  \eqref{odddiffeqn} we iterate \eqref{Gseq} over odd elements $g_{n-1}$ to $g_0$. We have that
\begin{equation*}
\begin{aligned}
g_{n} &= \dfrac{a}{2} g_{n-1} + \dfrac{b}{2} \\
&= \bigg(\frac{a}{2} \bigg) \bigg(\dfrac{a}{2}  g_{n-2} + \dfrac{b}{2} \bigg) + \dfrac{b}{2} \\
&= \bigg(\frac{a}{2} \bigg) \bigg(\dfrac{a}{2}  \Big(\dfrac{a}{2}  g_{n-3} + \dfrac{b}{2} \Big) + \dfrac{b}{2} \bigg) + \dfrac{b}{2} \\
&\vdotswithin{=} \\
& = \bigg(\frac{a}{2} \bigg)^n g_{0} + \dfrac{b}{2}
\Bigg( \bigg(\dfrac{a}{2}\bigg)^{n-1} + \bigg(\dfrac{a}{2}\bigg)^{n-2} + \ldots +
\bigg(\dfrac{a}{2}\bigg)^{2} + \bigg(\dfrac{a}{2}\bigg)^{1} + 1
 \Bigg) \\
&= \bigg(\frac{a}{2} \bigg)^n g_{0} + \dfrac{b}{2} \Bigg(
\dfrac{1-\big(a/2\big)^n}{1-\big(a/2\big)}
\Bigg) \quad \text{by the geometric series identity.}\\
\end{aligned}
\end{equation*}
The general solution to \eqref{odddiffeqn} for odd $g_0$ is thus
\begin{equation}\label{tempoddparity}
g_{n} = \bigg(\frac{a}{2} \bigg)^n \bigg( g_{0} + \frac{b}{a-2} \bigg) - \frac{b}{a-2}.
\end{equation}
Given that $g_0$ is an odd integer we now impose the constraint that $b \equiv 0$ mod ($a-2$) so that the factor $g_0+b/(a-2)$ in \eqref{tempoddparity} is an even integer. For odd $g_0$ we have that $g_n \equiv 0$ (mod $2$) when $n=\nu(g_0)$, yielding \eqref{oddparity}. 
\end{proof}
\end{thm}

\begin{remark}
Even integers $g_n$ in  \eqref{tempoddparity} may exist for odd $a$ and odd $b \not \equiv 0$ \textnormal{mod ($a-2$)}, but we do not investigate those sequences.
\end{remark}
The iteration exponents $\mu(g_0)$ and $\nu(g_0)$ in \eqref{evenparity} and \eqref{oddparity} are related to the $ruler$ function, which we now define.
\begin{defn}\label{rulerDefn}
For $n \in \mathbb{N}_{>0}$ the ruler function $r(n)$ is the exponent of the highest power of two that divides $2n$ (\seqnum{A001511} \cite{Sloane2018}) and satisfies the recursive equation
\begin{equation}\label{rulerFunction}
\begin{aligned}
r(2n) &= r(n)+1 \\
r(2n+1) &=1.
\end{aligned}
\end{equation}
\end{defn}
By Definition~\ref{rulerDefn} we have that $\mu(g_0) \in \mathbb{N}_{>0}$ and $\nu(g_0) \in \mathbb{N}_{>0}$ in \eqref{evenparity} and \eqref{oddparity}, respectively, satisfy 
\begin{equation}\label{xyEqn}
\begin{aligned}
\mu(g_0) &= r(g_0/2) && \text{if } g_0 \equiv \modd{0} {2} \\
\nu(g_0) &= r((g_0+b/(a-2))/2) && \text{if } g_0 \equiv \modd{1} {2}.
\end{aligned}
\end{equation}

We now present a generalized, alternating parity, Terras-modified Collatz function with the recursive ruler function $r(2n)=r(n)+1$ and $r(2n+1)=1$. We denote this function as the $GAPT$ Collatz function.
\begin{defn}[GAPT Collatz Function]
We let $c=b/(a-2)$ with odd $a \in \mathbb{N}$ and odd $b \in \mathbb{N}$ in \eqref{oddparity} such that $c$ is a positive odd integer.  For $n \in \mathbb{N}$, $x(n)  \in \mathbb{N}_{>0}$ \eqref{xyEqn}, $y(n) \in \mathbb{N}_{>0}$ \eqref{xyEqn}, and $h_n \in \mathbb{N}_{>0}$ we define the \textit{generalized alternating parity Terras-modified (GAPT) Collatz sequence} $H(h_0)= (h_0, h_1, \ldots)$ with the recursive function
\begin{equation}\label{Hseq}
h_{n+1} =
\begin{cases}
\bigg(\dfrac{1}{2}\bigg)^{r(h_n/2)} h_n & \text{if } h_n  \equiv \modd{0} {2} \\
\bigg(\dfrac{a}{2}\bigg)^{r((h_n+c)/2)} (h_n+c) - c &  \text{if } h_n \equiv \modd{1} {2}.
\end{cases}
\end{equation}
\end{defn}
As far as we know a similar ruler function representation for odd Collatz iterations first appears in the work of Crandall \cite{Cran1978}, however, an explicit functional form for $r(n)$ is not provided. We have found a similar (unpublished) approach \cite{Motta}.
\begin{thm}\label{finitenessThm}
Subsequences of consecutive even or odd elements in the generalized sequence $G$ \eqref{Gseq} are finite in length.
\begin{proof}
For generalized sequences $G$ \eqref{Gseq} and $H$ \eqref{Hseq} we have by construction that $H \lessdot G$. The exponents in the alternating parity sequence $H$ \eqref{Hseq} exist and are finite by Definition \ref{rulerDefn}. Thus, subsequences of consecutive even or odd elements in $G$ \eqref{Gseq} are finite in length.
\end{proof}
\end{thm}

\section{Alternating Parity Terras-modified (APT) Collatz sequences}

Henceforth we only consider the original $3n+1$ Collatz conjecture with $a=3$, $b=1$, and $c=1$ in \eqref{Hseq}. We now define the alternating parity Terras-modified (\textit{APT}) Collatz function.

\begin{defn}[APT Collatz Sequences]
For $n \in \mathbb{N}$ and $a_n \in \mathbb{N}_{>0}$ we define the \textit{alternating parity Terras-modified (APT) Collatz sequence} $A(a_0)= (a_0, a_1, \ldots)$ as
\begin{equation}\label{APTseq}
a_{n+1} =
\begin{cases}
\bigg(\dfrac{1}{2}\bigg)^{r(a_n/2)}a_n & \text{if } a_n  \equiv \modd{0} {2} \\
\bigg(\dfrac{3}{2}\bigg)^{r((a_n+1)/2)} (a_n+1) - 1 &  \text{if } a_n \equiv \modd{1} {2}
\end{cases}
\end{equation}
with $r(n)$ given by Definition \ref{rulerDefn}.
\end{defn}

\begin{coro}\label{finitenessCoro}
For $n \in \mathbb{N}_{>0}$ subsequences of consecutive even or odd elements in $C(n)$ \eqref{Cseq} and $T(n)$ \eqref{Tseq}  are finite in length.
\begin{proof}
Corollary~\ref{finitenessCoro} follows from Theorem~\ref{finitenessThm} with $a=3$, $b=1$, and $c=1$ in \eqref{Hseq}.
\end{proof}
\end{coro}

\begin{thm}\label{APTthms}
For $n \in \mathbb{N}_{>0}$ the sequences $C(n)$ \eqref{Cseq}, $T(n)$ \eqref{Tseq}, and $A(n)$ \eqref{APTseq} satisfy the following statements.
\begin{enumerate}[label=\alph*)]
\item $A(n) \lessdot T(n) \lessdot C(n)$.
\item $1 \in C(n)$ if and only if $1 \in A(n)$.
\item  If $1 \in C(n)$ then $|A_1(n)| \le |T_1(n)| \le |C_1(n)|$.
\end{enumerate}
\begin{proof}
These statements follow directly by construction of $A$ \eqref{APTseq} from $T$ \eqref{Tseq} and of $T$ \eqref{Tseq} from $C$ \eqref{Cseq} by reducing the number of iterations in $C$ and $T$.
\end{proof}
\end{thm}

\begin{conj}\label{mainConj}
For $n \in \mathbb{N}_{>0}$ the APT Collatz sequences $A(n)$ \eqref{APTseq} satisfy $1 \in A(n)$, implying that the \textit{CC} is true.
\end{conj}

\section{Modified APT (MAPT) Collatz sequences}

We now address the ruler function exponents in \eqref{APTseq}.  We first define two functions that are related to the ruler function to simplify (perhaps) this analysis.
\begin{defn}\label{grundyDef}
For $n \in \mathbb{N}$ we define the recursive function $p:\mathbb{N} \rightarrow \mathbb{N}$ as
\begin{equation}\label{pFunction}
\begin{aligned}
p(2n) &= n \\
p(2n+1) &= p(n)
\end{aligned}
\end{equation}
with $p(0) = 0$ (\seqnum{A0025480} \cite{Sloane2018}). Herein we also use the notation that $p_n=p(n)$.
\end{defn}
\begin{defn}\label{modRulerDefn}
For $n \in \mathbb{N}$ we define the recursive function $q:\mathbb{N} \rightarrow \mathbb{N}$ in terms of the ruler function $r(n)$ \eqref{rulerFunction} (\seqnum{A001511} \cite{Sloane2018})  as
\begin{equation*}
q(n) = r(n+1),
\end{equation*}
implying that
\begin{equation}\label{qFunction}
\begin{aligned}
q(2n) &= 1 \\
q(2n+1) &= q(n)+1
\end{aligned}
\end{equation}
with $q(0) = 1$. Herein we also use the notation that $q_n=q(n)$.
\end{defn}
We now represent even and odd integers in terms of $p_n$ \eqref{pFunction} and $q_n$ \eqref{qFunction}, although alternative representations are possible.
\begin{lem}\label{EvenOddLem}
For $n \in \mathbb{N}$, $p_n \in \mathbb{N}_{>0}$ and $q_n \in \mathbb{N}_{>0}$, the $n^{\text{\tiny th}}$ elements of the nondecreasing sequences of even positive integers $e_n=2(n+1)$ and odd positive integers $o_n=2n+1$ are
\begin{equation}\label{evenRep}
e_n (p_n, q_n) = (2p_n+1) 2^{q_n}
\end{equation}
and
\begin{equation}\label{oddRep}
o_n (p_n, q_n) = (2p_n+1) 2^{q_n} - 1,
\end{equation}
respectively.
\begin{proof}
For $m,n \in \mathbb{N}$ we prove by induction. For $n=0$ the base case for \eqref{evenRep} is
\begin{equation*}
\begin{aligned}
e_0 &= (2p_0+1)2^{q_0} \\
&= (2(0)+1)2^1 \\
&=2.
\end{aligned}
\end{equation*}
For the inductive step for \eqref{evenRep} we assume that there exists an $m$ such that $e_m = (2p_m+1) 2^{q_m} =2(m+1)$ and show that $e_{m+1} =2(m+2)$. For even $m$ we have odd $m+1=2j+1$ for some $j  = m/2 \in \mathbb{N}$ yielding
\begin{equation*}
\begin{aligned}
e_{m+1} &= (2p_{m+1}+1)2^{q_{m+1}} &\\
&= (2p_{2j+1}+1)2^{q_{2j+1}} &\\
&= (2p_{j}+1)2^{q_j+1} &\\
&= (2(j+1)) 2^1  &\quad \text{by the inductive hypothesis}\\
&= (2(m/2+1))2 &\\
&= 2(m+2).&\\
\end{aligned}
\end{equation*}
For odd $m$ we have even $m+1=2j$ for some $j  = (m+1)/2 \in \mathbb{N}$ yielding
\begin{equation*}
\begin{aligned}
e_{m+1} &= (2p_{m+1}+1)2^{q_{m+1}}\\
&= (2p_{2j}+1)2^{q_{2j}}\\
&= (2j+1)2^{1}\\
&= (2((m+1)/2)+1)2^{1}\\
&= 2(m+2).
\end{aligned}
\end{equation*}
For $n = 0$ the base case for \eqref{oddRep} is
\begin{equation*}
\begin{aligned}
o_0 &= (2p_0+1)2^{q_0} - 1\\
&= (2(0)+1)2^1 - 1\\
&=1.
\end{aligned}
\end{equation*}
For the inductive step for \eqref{oddRep} we assume that there exists an $m$ such that $o_m = (2p_m+1) 2^{q_m}-1 =2m+1$ and show that $o_{m+1} =2m+3$. For even $m$ we have odd $m+1=2j+1$ for some $j  = m/2 \in \mathbb{N}$ yielding
\begin{equation*}
\begin{aligned}
o_{m+1} &= (2p_{m+1}+1)2^{q_{m+1}} - 1 & \\
&= (2p_{2j+1}+1)2^{q_{2j+1}} - 1 & \\
&= (2p_{j}+1)2^{q_j+1} - 1 & \\
&= (2(j+1)) 2^1 - 1 &\quad \text{by the inductive hypothesis} \\
&= (m+2) 2^1 - 1 & \\
&= 2m+3. & \\
\end{aligned}
\end{equation*}
For odd $m$ we have even $m+1=2j$ for some $j  = (m+1)/2 \in \mathbb{N}$ yielding
\begin{equation*}
\begin{aligned}
o_{m+1} &= (2p_{m+1}+1)2^{q_{m+1}} - 1 \\
&= (2p_{2j}+1)2^{q_{2j}} - 1\\
&= (2j+1)2^{1} - 1 \\
&= (m+2)2^1 - 1 \\
&= 2m+3.
\end{aligned}
\end{equation*}
\end{proof}
\end{lem}

We substitute $e_n$ \eqref{evenRep} and $o_n$ \eqref{oddRep} for even and odd $a_n$, respectively, into the ruler function arguments in \eqref{APTseq} and obtain the following identities.
\begin{lem}\label{rulerLemma}
For even and odd integers $e_n$ \eqref{evenRep} and $o_n$ \eqref{oddRep}, respectively, the ruler function \eqref{rulerFunction} satisfies the identities
\begin{equation}\label{rulerIdentEven}
r(e_n (p_n, q_n)/2) = q_n
\end{equation}
and
\begin{equation}\label{rulerIdentOdd}
r((o_n (p_n, q_n)+1)/2) = q_n.
\end{equation}
\begin{proof}
For $n \in \mathbb{N}$ and $q_n \ge 1$ we recurse the argument of the ruler function \eqref{rulerFunction} and have that
\begin{equation*}
\begin{aligned}
r(e_n (p_n,q_n)/2) &= r((2p_n+1) 2^{{q_n}-1}) \\
&= r((2p_n+1) 2^{q_n-2}) + 1 \\
&= r((2p_n+1) 2^{q_n-3}) + 2 \\
& \vdotswithin{=} \\
&= r((2p_n+1) 2^{q_n-q_n}) + q_n-1 \\
&= r(2p_n+1) + q_n-1 \\
&= q_n \\
\end{aligned}
\end{equation*}
and, thus, that
\begin{equation*}
\begin{aligned}
r((o_n (p_n, q_n)+1)/2) &= r((2p_n+1) 2^{q_n-1}) \\
&=q_m. \\
\end{aligned}
\end{equation*}
\end{proof}
\end{lem}

Lemmas~\ref{EvenOddLem} and \ref{rulerLemma} applied to the APT Collatz sequences  \eqref{APTseq} yield a modified representation of \eqref{APTseq}. We denote this modification as the $MAPT$ Collatz function.
\begin{thm}[MAPT Collatz Function] \label{APTmap}
For $i,j,n \in \mathbb{N}$ and $p_i$ and $q_i$ in Definitions~\ref{grundyDef} and \ref{modRulerDefn}, respectively, the APT Collatz function \eqref{APTseq} is equivalent to the maps of an even $a_n(i)$ to an odd $a_{n+1}(i)$
\begin{equation}\label{aniseven}
\begin{aligned}
a_n (i) &= (2p_i+1) 2^{q_i}, & (a_n &\equiv \modd{0} {2})\\
a_{n+1}(i) &= (2 p_i + 1), & (a_{n+1} &\equiv \modd{1} {2})
\end{aligned}
\end{equation}
and an odd $a_{n}(j)$ to an even $a_{n+1}(j)$
\begin{equation}\label{anisodd}
\begin{aligned}
a_{n}(j) &= (2p_j+1) 2^{q_j} - 1, & (a_n &\equiv \modd{1} {2} ) \\
a_{n+1}(j) &= (2p_j+1) 3^{q_j} -1, & (a_{n+1} &\equiv \modd{0} {2}).
\end{aligned}
\end{equation}
\begin{proof}
For $n, a_n \in \mathbb{N}$ and for some $i,j \in \mathbb{N}$ in \eqref{APTseq} we let even $a_n(i) = e_i (p_i, q_i)$ \eqref{evenRep} and odd $a_n(j) = o_i (p_j, q_j)$ \eqref{oddRep} from Lemma~\ref{EvenOddLem}, and we use the ruler function identities \eqref{rulerIdentEven} and \eqref{rulerIdentOdd} in Lemma~\ref{rulerLemma}.  For even $a_n(i)$ we have odd $a_{n+1}(i)$
\begin{equation*}
\begin{aligned}
a_{n+1}(i) &= \bigg(\dfrac{1}{2}\bigg)^{r(a_n(i)/2)} a_n(i) \\
&= \bigg(\dfrac{1}{2}\bigg)^{r(e_i (p_i,q_i)/2)} e_i (p_i,q_i) \\
&= \bigg(\dfrac{1}{2}\bigg)^{q_i}  (2p_i+1)  2^{q_i} \\
&= 2p_i+1.
\end{aligned}
\end{equation*}
For odd $a_n(j)$ we have even $a_{n+1}(j)$
\begin{equation*}
\begin{aligned}
a_{n+1}(j) &= \bigg(\dfrac{3}{2}\bigg)^{r((a_n(j)+1)/2)} (a_n(j)+1) - 1  \\
&=\bigg(\dfrac{3}{2}\bigg)^{r((o_j (p_j, q_j)+1)/2)} (o_j (p_j, q_j)+1) - 1  \\
&=\bigg(\dfrac{3}{2}\bigg)^{q_j} ((2p_j+1) 2^{q_j} - 1+1) - 1  \\
&=(2p_j+1) 3^{q_j} -1.
\end{aligned}
\end{equation*}
\end{proof}
\end{thm}

\section{Even (EMAPT) and Odd (OMAPT) MAPT Collatz sequences}

We now generate even-only and odd-only $MAPT$ Collatz subsequences of $A$ \eqref{APTseq} as represented in Equations~\eqref{aniseven} and \eqref{anisodd} in Theorem~\ref{APTmap}, respectively. 

\begin{thm}\label{EvenOddThm}
The subsequences of even elements of $A$ \eqref{APTseq} are
\begin{equation}\label{evenA}
a_{n+2} = \bigg(2  p\bigg(p\big((a_{n}-2)/2\big)\bigg)+1\bigg)  3^{q\big(p\big((a_{n}-2)/2\big)\big)} -1 \\
\end{equation}
with $a_{n} = 2m$ for $m \in \mathbb{N}$, and the subsequences of odd elements of $A$ \eqref{APTseq} are
\begin{equation}\label{oddA}
a_{n+3} = 2  p \bigg( \bigg( \Big( 2 p\big((a_{n+1}-1)/2\big)+1 \Big)3^{q\big((a_{n+1}-1)/2\big)} - 3 \bigg)/2 \bigg)+1 \\
\end{equation}
with $a_{n+1} = 2m+1$ for $m \in \mathbb{N}$.
\begin{proof}
We prove algebraically. We recurse Equations~\eqref{aniseven} and \eqref{anisodd} in Theorem~\ref{APTmap} and use the results of Lemma~\ref{EvenOddLem}. For $i,j,k, \ldots \in \mathbb{N}$ we start with some $i$ and even $a_n(i)$ yielding
\begin{equation*}
\begin{aligned}
a_n(i) &= (2p_i+1) 2^{q_i} \\
&=2(i+1) \\
a_{n+1}(i) &= 2p_i+1   \\
&=a_{n+1}(j)\\
a_{n+1}(j)  &= (2p_j+1)2^{q_j} - 1 \\
&= 2j+1 \\
a_{n+2}(j)  &= (2p_j+1)3^{q_j} - 1 \\
&= a_{n+2}(k) \\
a_{n+2}(k) &= (2p_k+1) 2^{q_k} \\
&= 2(k+1) \\
a_{n+3}(k) &= 2p_k+1 \\
&= \ldots.
\end{aligned}
\end{equation*}
For even $a_n(i)$ we have that $i=(a_n-2)/2$ and $j=p_i$. We substitute these values for $i$ and $j$ into $a_{n+2}(j)= (2p_j+1)3^{q_j} - 1$ yielding \eqref{evenA}. For odd $a_{n+1}(j)$ we have that $j=(a_{n+1}-1)/2$ and $k=((2p_j+1)3^{q_j} - 3)/2$. We substitute these values for $j$ and $k$ into $a_{n+3}(k)=2p_k+1$ yielding \eqref{oddA}.
\end{proof}
\end{thm}
For even $n \in \mathbb{N}$ we define the even-only elements of sequences $A(n)$ by sequences $U(n)$.
\begin{defn} [EMAPT Collatz Function]\label{EMAPTdefn}
By Theorem~\ref{EvenOddThm} we denote the sequences of even positive integers of $A$ \eqref{evenA} by $U$ (Even MAPT Collatz function) with even $u_n=a_n$ and even $u_{n+1}=a_{n+2}$ yielding
\begin{equation}\label{Useq}
u_{n+1} = \big(2 p\big(p\big((u_{n}-2)/2\big)\big)+1\big) 3^{q\big(p\big((u_{n}-2)/2\big)\big)} -1 \\
\end{equation}
with $u_0=2m$ for $m \in \mathbb{N}$.
\end{defn}

For odd $n \in \mathbb{N}$ we define the odd-only elements of $A(n)$ by sequences $V(n)$.
\begin{defn} [OMAPT Collatz Function]
By Theorem~\ref{EvenOddThm} we denote the sequences of odd positive integers of $A$ \eqref{oddA} by $V$ (Odd MAPT Collatz function) with odd $v_n=a_{n+1}$ and odd $v_{n+1}=a_{n+3}$ yielding
\begin{equation}\label{Vseq}
v_{n+1} = 2 p \big( \big( \big( 2 p\big((v_{n}-1)/2\big)+1 \big)3^{q\big((v_{n}-1)/2\big)} - 3 \big)/2 \big)+1 \\
\end{equation}
with $v_0=2m+1$ for $m \in \mathbb{N}$.
\end{defn}
For even $u_n$ and odd $v_n$ elements of sequences $U$~\eqref{Useq} and $V$~\eqref{Vseq}, respectively, we have that
\begin{equation}\label{UVrelation}
v_n = \frac{2u_n}{2^{r(u_n)}}
\end{equation}
by $A$ \eqref{APTseq}.  The sequences $U$ and $V$ satisfy the covering relations $U(n) \lessdot A(n)$ and $V(n)  \lessdot A(n)$ for even and odd $n \in \mathbb{N}$, respectively.  For even $n \in \mathbb{N}_{>0}$ if sequences $U(n)$ \eqref{Useq} satisfy $2 \in U(n)$, then the \textit{CC} is true.  Equivalently, by \eqref{UVrelation} and for odd $n \in \mathbb{N}_{>0}$, if sequences $V(n)$ \eqref{Vseq} satisfy $1 \in V(n)$, then the \textit{CC} is true.

\section{EMAPT Collatz sequences: $U(6n+2)$}

We now constrain the set of even integers $u_0$ in sequences $U(u_0)$ \eqref{Useq} sufficient to prove the \textit{CC}.
\begin{thm}\label{6m2Thm}
For $n \in \mathbb{N}$ and even $u_n \in \mathbb{N}_{>0}$ Equation~\eqref{Useq} has the form $u_{n+1}= 6m+2$ for some $m \in \mathbb{N}$.
\begin{proof}
For $n \in \mathbb{N}$ we define the functions
$f_1(n) \in \mathbb{N}_{>0}$ and $f_2(n) \in \mathbb{N}_{>0}$ as
\begin{equation*}
f_1(n) = p\big(p\big((u_{n}-2)/2\big)\big)
\end{equation*}
and
\begin{equation*}
f_2(n) = q\big(p\big((u_{n}-2)/2\big)\big)-1.
\end{equation*}
For some $m \in \mathbb{N}$ Equation~\eqref{Useq} has the form
\begin{equation*}
\begin{aligned}
u_{n+1} &= (2 f_1(u_n)+1) 3^{f_2(u_n)+1} -1 \\
&= 6\left(3^{f_2(u_n)} f_1(u_n) + \dfrac{3^{f_2(u_n)}-1}{2}\right)+2 \\
&= 6m+2.
\end{aligned}
\end{equation*}
\end{proof}
\end{thm}

By Theorem~\ref{6m2Thm} we only need to show that $2 \in U$ for even numbers of the form $u_0=6m+2$ to satisfy the \textit{CC}.  A similar result appears in the unpublished work \cite{Motta}. We state this result.

\begin{thm}\label{UThm}
For $n \in \mathbb{N}$ if the EMAPT Collatz sequences $U(6n+2)$ \eqref{Useq} satisfy $2 \in U$ then the \textit{CC} is true.
\begin{proof}
This theorem follows from Theorem~\ref{6m2Thm}.
\end{proof}
\end{thm}

\section{The EMAPT Collatz sequences beyond $U(6n+2)$}

For $n \in \mathbb{N}$, $\alpha \in \mathbb{N}$, and $\beta \in \mathbb{N}$ we now further constrain the set of the integer arguments of $U(\alpha n + \beta)$ \eqref{Useq} sufficient to prove the \textit{CC}.  Indeed, Monks (2006) \cite{Monks2006} shows that the set of arguments $\{\alpha n + \beta \}$ is sufficient to prove the \textit{CC} for any nonnegative integers $\alpha$ and $\beta$ with $\alpha \ne 0$.

For $n \in \mathbb{N}$ we consider the elements of $U$ \eqref{Useq} and reproduce here as
\begin{equation*}
u_{n+1} = \big(2 p\big(p\big((u_{n}-2)/2\big)\big)+1\big) 3^{q\big(p\big((u_{n}-2)/2\big)\big)} -1.
\end{equation*}
For index $x_n \in \mathbb{N}$ we let $u_n=6x_n+2$ and $u_{n+1}=6x_{n+1}+2$. Theorem~\ref{6m2Thm} implies that
\begin{equation*}
6x_{n+1}+2= \Big(2 \, p\big(p\big(3x_n\big)\big)+1\Big) \, 3^{q\big(p\big(3x_n\big)\big)} -1,
\end{equation*}
and yields the sequence $X$ of indices of $U$
\begin{equation}\label{Xseq}
x_{n+1} = \frac{1}{2} \left( \Big(2 \, p\big(p\big(3x_n\big)\big)+1\bigg) \, 3^{q\big(p\big(3x_n\big)\big)-1} -1\right).
\end{equation}
We may thus state the following theorem.
\begin{thm}\label{XThm}
For $n \in \mathbb{N}$ and $x_n \in \mathbb{N}$ in \eqref{Xseq}  if sequences $X$ satisfy $0 \in X$ then the \textit{CC} is true.
\begin{proof}
For $n \in \mathbb{N}$ and $x_n \in \mathbb{N}$ if $x_n = 0$, then $2 \in U$.
\end{proof}
\end{thm}
We can strengthen Theorems~\ref{UThm} and \ref{XThm}. For $n \in \mathbb{N}$, $x_n \in \mathbb{N}$, and $m \in \mathbb{N}$, we observe computationally that $x_{n+1} \ne 3m+2$.  Before proving that statement, we need the following lemma.

\begin{lem}\label{p3nparity}
The function $p:\mathbb{N} \rightarrow \mathbb{N}$ with $p(2n)=n$ and $p(2n+1)=p(n)$ \eqref{pFunction} satisfies $p(3n) \equiv 0$ \textnormal{(mod $3$)} or $p(3n) \equiv 2$ \textnormal{(mod $3$)} for all $n \in \mathbb{N}$.
\begin{proof}
We use strong induction. For the base case we have that $p(0)=0$. For the inductive step we assume for some arbitrary $n \in \mathbb{N}$ and for $m \in \mathbb{N}$ that $p(3m) \equiv 0$ (mod $3$) or that $p(3m) \equiv 2$ (mod $3$) for $m \in \{0,1,2,\ldots,n\}$.

An integer $n \in \mathbb{N}$ is even or odd.  Thus, an integer $a \in \mathbb{N}$ exists such that $n=2a$ or $n=2a+1$.  For even $n=2a$ we have that
\begin{equation*}
p(3(n+1)) =  p(3(2a+1)) = p(6a+3)=p(2(3a+1)+1)=p(3a+1).
\end{equation*}
Similarly, an integer $b \in \mathbb{N}$ exists such that $a=2b$ or $a=2b+1$. For even $a=2b$ we have that
\begin{equation*}
p(3(n+1)) = p(3(2b)+1) =  p(3b) \equiv 0 \,\, \text{or } 2 \pmod{3}
\end{equation*}
since $n=2a=2(2b)$ and $b=n/4 \in \{0,1,\ldots,n\}$.
For odd $a=2b+1$ we have that
\begin{equation*}
p(3(n+1)) = p(3(2b+1)+1) =  p(2(3b+2)) = 3b+2 \equiv 2 \pmod{3}.
\end{equation*}
For odd $n=2a+1$ we have that
\begin{equation*}
p(3(n+1)) =  p(3((2a+1) +1)) = p(6a+6)=p(2(3a+3))=3a+3 \equiv 0 \pmod{3}.
\end{equation*}
\end{proof}
\end{lem}

We now proceed constraining indices $x_n$ in sequence $X$ \eqref{Xseq}.
\begin{thm}\label{XnThm}
For $x_n \in \mathbb{N}$ elements $x_{n+1}$ in \eqref{Xseq} satisfy $x_{n+1} \ne 3m+2$ for $m \in \mathbb{N}$; equivalently, $x_{n+1} \equiv 0$ \textnormal{(mod } $3$\textnormal{)} or  $x_{n+1} \equiv 1$ \textnormal{(mod }$3$\textnormal{)}.
\begin{proof}
We prove by contradiction.  For $x_n \in \mathbb{N}$ we assume that $x_{n+1}=3m+2$ for some $m \in \mathbb{N}$ in \eqref{Xseq} yielding
\begin{equation*}
3m+2 = \frac{1}{2} \left( \Big(2 \, p\big(p\big(3x_n\big)\big)+1\Big) \, 3^{q\big(p\big(3x_n\big)\big)-1} -1\right)
\end{equation*}
and, thus,
\begin{equation*}
6m+5 = \Big(2 \, p\big(p\big(3x_n\big)\big)+1\Big) \, 3^{q\big(p\big(3x_n\big)\big)-1}.
\end{equation*}
By Equations~\eqref{pFunction} and \eqref{qFunction} we have that $q(p(3x_n)) \ge 1$. If $q(p(3x_n)) > 1$ then we have a contradiction because $3 \nmid (6m+5)$. If $q(p(3x_n))=1$ then we have even $p(3x_n)$ and $p(p(3x_n))=p(3x_n)/2$. We thus have that
\begin{equation*}
6m+4 = p(3x_n).
\end{equation*}
By Lemma~\ref{p3nparity} for some $i \in \mathbb{N}$ we have that $p(3x_n)=3i$ or $p(3x_n)=3i+2$, implying that $6m+4 =3i$ or $6m+4=3i+2$, respectively. We have contradictions in both cases because $3 \nmid (6m+4)$ and $3 \nmid (6m+2)$, respectively.
\end{proof}
\end{thm}
We can thus state the one of the main results of this work.
\begin{thm}\label{18mUThm}
For $u_0=18n+2$ or $u_0=18n+8$ in $U$ \eqref{Useq} if $2 \in U$, then the \textit{CC} is true.
\begin{proof}
The proof follows directly from Theorem~\ref{XnThm}.
\end{proof}
\end{thm}

\section{The EMAPT Collatz sequence $U(n)$ in terms of $r(n)$ }

For $m \in \mathbb{N}$ we now express the recursive terms $p(p(m))$ and $q(p(m))$ in \eqref{Useq} as a function of the ruler function $r(m)$ \eqref{rulerFunction}.  For $i \in \mathbb{N}$ and $j \in \mathbb{N}$ we uniquely represent the integers $m$ in terms of the two-tuples $(i,j)$ as
\begin{equation*}\label{ijRep}
m = (2i+1) 2^j - 1.
\end{equation*}
For $k \in \mathbb{N}$ and $l \in \mathbb{N}$ we represent the integers $i$ in terms of the two-tuples $(k,l)$ as
\begin{equation*}\label{klRep}
i = (2k+1) 2^l - 1.
\end{equation*}
Thus, we have the three-tuple $(j,k,l)$ representation of $m$
\begin{equation*}\label{ijklRep}
m = (2 ( (2k+1) 2^l - 1 ) + 1) 2^j - 1.
\end{equation*}
We now recurse $p(p(m))$ and $q(p(m))$ using the definitions for $p(m)$ \eqref{pFunction} and $q(m)$ \eqref{qFunction}. We obtain that
\begin{equation*}
\begin{aligned}
p(p(m)) &= k \\
q(p(m)) &= l+1.
\end{aligned}
\end{equation*}
Thus, we have for any integer $(u_n-2)/2$ that
\begin{equation*}
(u_n(j,k,l) - 2)/2 = \left(2( (2k+1)2^l - 1) + 1\right)2^j - 1
\end{equation*}
and, thus, that
\begin{equation*}
u_n(j,k,l) = \left( (2k+1)2^{l+1} - 1\right)2^{j+1}.
\end{equation*}
By \eqref{Useq} we have that
\begin{equation*}
u_{n+1}(k,l) = (2k +1) \, 3^{l+1} -1.
\end{equation*}
We eliminate the common factor $2k+1$ and have that
\begin{equation}\label{jlEqn}
u_{n+1}(j,l) = \left(\dfrac{u_n}{2^{j+1}}+1\right)\left(\dfrac{3}{2}\right)^{l+1}-1.
\end{equation}
The \textit{maximum allowed} tuples $(j,l)$ are the $2$-adic valuations of factors in \eqref{jlEqn}. We have that
\begin{equation*}
\begin{aligned}
j &= r(u_n) -2 \\
l &= r\left(\dfrac{u_n}{2^{r(u_n)-1}}+1\right)-2.
\end{aligned}
\end{equation*}
Thus, a primary result of this work is that
\begin{equation}\label{foobar}
u_{n+1} =  \left(\dfrac{2u_n}{2^{r(u_n)}}+1\right)\left(\dfrac{3}{2}\right)^{r\left(\dfrac{2u_n}{2^{r(u_n)}}+1\right)-1}-1
\end{equation}
with \eqref{UVrelation}
\begin{equation*}
v_n = \frac{2u_n}{2^{r(u_n)}}.
\end{equation*}
For even $u_0 \in \mathbb{N}_{>0}$ in $U(u_0)$ \eqref{foobar} all subsequent elements $u_n$ in sequence $U(u_0)$ are of the form $18m+2$ or $18m+8$. (We note that for odd $u_0 \in \mathbb{N}_{>0}$ in $U(u_0)$ \eqref{foobar} all subsequent elements $u_n$ beyond $u_1$ in sequence $U(u_0)$ appear to be of the form $18m+2$ or $18m+8$.) We conjecture that $2 \in U(u_0)$, which would imply that the $CC$ is true.

If we substitute $v_n$ \eqref{UVrelation} into \eqref{foobar}, then we more clearly have a fixed rational root for the relation between adjacent elements in the sequences $U$ and $V$
\begin{equation}\label{fixedPt}
\left(\dfrac{u_{n+1}+1}{v_n+1}\right)^{1/r((v_n+1)/2)} = \dfrac{3}{2}.
\end{equation}

\section{$Z$: Reverse Map of $U$}

For $n \in \mathbb{N}$ and $\alpha_n, \beta_n \in \mathbb{N}_{>0}$ we let $\alpha_n =r(v_n+1)-1$ and $\beta_n = r(u_n)-1$. The infinite linear system of equations \eqref{foobar} for adjacent elements of $U$ satisfies
\begin{equation}\label{linsysU}
2^{\alpha_n+\beta_n} u_{n+1} = 3^{\alpha_n} u_n + 2^{\beta_n}(3^{\alpha_n}-2^{\alpha_n}). 
\end{equation}
We now consider the linear system of equations for the reverse map $u_{n+1} \mapsto u_n$ in \eqref{linsysU}
\begin{equation}\label{revmapU}
u_n = 2^{\alpha_n+\beta_n} 3^{-\alpha_n} u_{n+1} - 2^{\beta_n}3^{-\alpha_n}(3^{\alpha_n}-2^{\alpha_n}).
\end{equation}
In \eqref{revmapU} we map $u_n \mapsto z_{n+1}$ and $u_{n+1} \mapsto z_n$ yielding
\begin{equation}\label{znrecursion}
z_{n+1} = m_{n+1} z_n + b_{n+1}
\end{equation}
with rational slopes
\begin{equation*}\label{mslope}
m_n = 2^{\alpha_n+\beta_n} 3^{-\alpha_n}
\end{equation*}
and rational intercepts
\begin{equation*}\label{bintercept}
b_n = - 2^{\beta_n}3^{-\alpha_n}(3^{\alpha_n}-2^{\alpha_n}).
\end{equation*}
Analogous to the $U$ and $V$ sequences, we define the odd only $W$ sequences $w_n \ne 0 \pmod{3}$ by
\begin{equation}\label{wneqn}
w_n = \frac{2z_n}{2^{r(z_n)}}
\end{equation}
and
\begin{equation}\label{zneqn}
\dfrac{z_{n}+1}{w_{n+1}+1} = \left(\dfrac{3}{2}\right)^{r((w_{n+1}+1)/2)}
\end{equation}
with $\alpha_n =r(w_{n+1}+1)-1 = r((w_{n+1}+1)/2)$ and $\beta_n = r(z_n)-1$.
\begin{thm}\label{z0thm}
($z_0$-cycle; Steiner \cite{Steiner1977}) If $z_1 = z_0$ in \eqref{znequation} then $z_0=2$.
\begin{proof}
In the original $CC$ sequence the only 1-cycle is the trivial sequence $(2,1,2,1, \ldots)$ \cite{Steiner1977}. We assume that $z_1 = z_0$. We have by \eqref{znrecursion} that
\begin{equation}\label{zcycle}
\begin{aligned}
z_0 &= \dfrac{b_1}{1-m_1} \\
&=\dfrac{- 2^{\beta_1}3^{-\alpha_1}(3^{\alpha_1}-2^{\alpha_1})}{1-2^{\alpha_1+\beta_1} 3^{-\alpha_1}} \\
&=\dfrac{- 2^{\beta_1}(3^{\alpha_1}-2^{\alpha_1})}{3^{\alpha_1}-2^{\alpha_1+\beta_1}}. \\
\end{aligned}
\end{equation}
The only 2-tuples that satisfy \eqref{zcycle} with $z_0 \in \mathbb{N}$ are $(\alpha_1,\beta_1)=(1,1)$, and we have that $z_0=2$.
\end{proof}
\end{thm}
We can select an initial $z_0$ and recurse \eqref{znrecursion} yielding
\begin{equation}\label{znequation}
z_n = \left(
\displaystyle\prod_{i=1}^{n} m_i \right) z_0 + b_{n} + \displaystyle\sum_{i=1}^{n-1} \left( b_i \displaystyle\prod_{j=i+1}^{n} m_j \right).
\end{equation}
We can let
\begin{equation*}\label{Mequation}
M_n (z_n) = \left(
\displaystyle\prod_{i=1}^{n} m_i \right)
\end{equation*}
and
\begin{equation*}\label{Bequation}
B_n (z_n) = b_{n} + \displaystyle\sum_{i=1}^{n-1} \left( b_i \displaystyle\prod_{j=i+1}^{n} m_j \right).
\end{equation*}
We thus have that
\begin{equation}\label{ZMB}
z_n = M_n(z_n) z_0 + B_n(z_n)
\end{equation}
with real 2-tuples $(M_n(z_n),B_n(z_n))$ for an infinite number of $z_n$ and possibly unique $w_n$ with $z_0=2$.  Rather than attempting to show that \eqref{ZMB} has $complete$ solutions in the evens with $z_0=2$, we now focus on \eqref{zneqn} to prove the Collatz conjecture.

\section{The $WZ$ Tree}

\begin{figure}[H]
\includegraphics[width=0.70\textwidth]{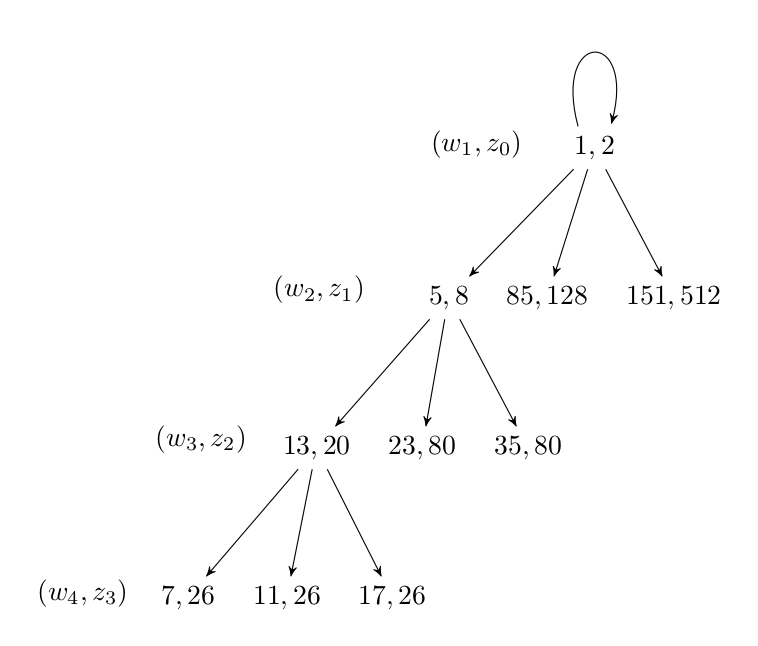}
\centering
\caption{Partial (upper left) WZ Tree}
\centering
\label{fig:WZtree}
\end{figure}

Herein we assume a fundamental knowledge of graphs and trees. We first illustrate some properties of the $WZ$ sequence tree, partially shown in Fig.~\ref{fig:WZtree} with root vertex $(w_1,z_0)=(1,2)$. Each 2-tuple $(w_n,z_n)$ connected by a directed edge satisfies \eqref{wneqn}. Each 2-tuple vertex $(w_{n+1},z_{n})$ satisfies \eqref{zneqn}. Each vertex has an infinite number of children; clearly, but see \eqref{infwnlem} below. Each level $n$ is ordered from left to right in terms of increasing unique $w_{n}$ terms and extends infinitely to the right; clearly again, but see \eqref{infwnlem} below.

For $m \in \mathbb{N}_{>0}$ the ordered sequence of all possible values for $w_{n+1}$ in \eqref{zneqn} is $(w_{n+1}(m)) = ((6m + (-1)^{m} - 3)/2)$ and is complete in the odds not divisible by $3$ (\seqnum{A007310} \cite{Sloane2018}). For each element in $(w_{n+1})$ we generate a $z_n$ \eqref{zneqn}, a $w_{n}$ \eqref{wneqn}, and the map $w_{n} \mapsto w_{n+1}$
\begin{equation}\label{wsets}
\begin{aligned}
(w_{n+1}) = (&1,5,7,11,13,17,19,23,25,29,31,35,37,41,43,47,49,53,55,59,\ldots) \\
(w_{n}) = (&1,1,13,13,5,13,11,5,19,11,121,5,7,31,49,121,37,5,47,67,\ldots). \\
\end{aligned}
\end{equation}

We now state two Lemmas.

\begin{lem}\label{cyclelem}
The only cycle in the Collatz tree is the trivial $1 \mapsto 1$.
\begin{proof}
The $w_{n+1}$ values are unique in $(w_{n+1})$ \eqref{wsets}, indicating that each vertex generated in the forward (down) direction in the Fig.~\ref{fig:WZtree} is unique except for the trivial cycle $1 \mapsto 1$ (See Theorem~\ref{z0thm}). 
\end{proof}
\end{lem}
\noindent
Given a $w_{n+1}$ the associated $w_n$ values are not uniquely determined, cover the set of integers $\{w_{n+1}\}$, and occur an infinite number of times.
\begin{lem}\label{infwnlem}
Each integer in the sequence of unique odds $(w_{n+1})$ occurs an infinite number of times in the sequence $(w_{n})$. 
\begin{proof}
For $m,n \in \mathbb{N}$ and a chosen $(z_n,w_n)$ 2-tuple, if $z_m \ne z_n$ and $w_m = w_n$ then, by \eqref{wneqn}, we have that
\begin{equation*}
\frac{2z_m}{2^{r(z_m)}} = \frac{2z_n}{2^{r(z_n)}}
\end{equation*}
and thus that
\begin{equation*}
z_m = 2^k z_n
\end{equation*}
for $k \in \mathbb{N}$ satisfying
\begin{equation*}
k = 2^{r(z_m)-r(z_n)}.
\end{equation*}
The sequence $(r(z_m))$ for all $z_m$ contains all positive integers an infinite number of times by Definition \ref{rulerDefn}. Thus, we have an infinite number of $w_m = w_n$.
\end{proof}
\end{lem}

\section{A Collatz Conjecture Proof}

We now present our last result.
\begin{thm}\label{ccproof}
The Collatz Conjecture is true.
\begin{proof}
We have a well--ordered, countable infinite set of unique $w_{n+1}$ vertices \eqref{wsets} that is complete in the odds not divisible by $3$. Except for the trivial cycle we also have a countable infinite set of unique $(w_{n+1},w_{n})$ acyclic edges (\eqref{cyclelem} and \eqref{infwnlem}). By Zorn's lemma, the well--ordering theorem or the axiom of choice (which are equivalent \cite{Jech1973}) we can construct an acyclic tree that spans \cite{Serre2003} our desired vertices, i.e., all odd numbers not divisible by $3$ with root vertex $1$. The existence of this tree with the covering relation $W(n)\lessdot C(n)$ is sufficient to complete the proof.
\end{proof}
\end{thm}

\section{Acknowledgements}

Thank you.

\end{document}